\documentclass[11pt,letterpaper]{article}
\usepackage{amsthm}
\usepackage{amsmath}
\usepackage{amsfonts}
\usepackage{latexsym}
\usepackage{geometry}
\usepackage{enumerate}

\emergencystretch=\maxdimen
\title{Entropy, noncollapsing, and a gap theorem for ancient solutions to the Ricci flow}
\author{Yongjia Zhang}
\numberwithin{equation}{section}
\begin{document}
\maketitle

In this paper we discuss the asymptotic entropy for ancient solutions to the Ricci flow. We prove a gap theorem for ancient solutions, which could be regarded as an entropy counterpart of Yokota's work. In addition, we prove that under some assumptions on one time slice of a complete ancient solution with nonnegative curvature operator, finite asymptotic entropy implies $\kappa$-noncollapsing on all scales. This provides an evidence for Perelman's more general assertion that on a complete ancient solution with nonnegative curvature operator, bounded entropy is equivalent to $\kappa$-noncollapsing.

\tableofcontents

\section{Introduction}
\newtheorem{theorem_gap}{Theorem}[section]
\newtheorem{theorem_noncollapsing}[theorem_gap]{Theorem}

A solution $(M^n,g(t))$ to the Ricci flow
\begin{displaymath}
\frac{\partial}{\partial t}g(t)=-2Ric(g(t))
\end{displaymath}
is called \textit{ancient} if $g(t)$ is defined on the time interval $(-\infty,T]$, for some $T>-\infty$. It is convenient to let $T=0$. Ancient solutions arise naturally as blow up limits of finite time singular solutions of the Ricci flow, and therefore are of great interest for the study of singularity formation. By analyzing ancient solutions with nonnegative curvature operator, Perelman \cite{perelman2002entropy} showed how singularities form in the Ricci flow of dimension three, which paved Perelman's path of solving the geometrization conjuecture via Hamilton's program \cite{hamilton1995formation}. However, general ancient solutions in higher dimensions are yet little understood. There are some fundamental results such as \cite{naber2010noncompact}, \cite{yokota2009perelman}, both of which work with reduced geometry. As to the role played by entropy in Type I ancient solutions, please refer to \cite{cao2011conjugate}, \cite{xuequation}.
\\

The entropy formula for the Ricci flow was introduced by Perelman \cite{perelman2002entropy}. He used it to prove his noncollapsing and pseudolocality theorems. A Ricci flow $(M^n,g(t))$ is called (weakly) $\kappa$-noncollapsed at scale $r_0$, if for all $(x,t)$, $r\leq r_0$, an upper bound of curvature  $|Rm|\leq r^{-2}$ on the parabolic cube $B_t(x,r)\times [t-r^2,t]$ implies a lower bound of the volume of the ball $\operatorname{Vol}(B_t(x,r))\geq \kappa r^n$. In this paper, we will always use $B_t(x,r)$ to denote the ball of radius $r$, centered at $x$, and with respect to the metric $g(t)$. It is well known that for a parabolic cube with a curvature bound as described above, a lower bound of the reduced volume centered at $(x,t)$ evaluated at $t-r^2$ implies a lower volume bound as given above, where the $\kappa$ depends only on the reduced volume and the dimension; see \cite{perelman2002entropy}. In this paper, we show that the same holds for Perelman's pointed entropy. There is also a notion of strong $\kappa$-noncollapsing, which assumes only a scalar curvature upper bound on the ball $B_t(x,r)$ (instead of the parabolic cube). In this paper, every noncollapsing notion we will mention means weak noncollapsing. Moreover, when dealing with ancient solutions, we always work with noncollapsing on all scales.
\\

Let $\bar{W}(x,t)$ denote the \textit{asymptotic entropy}, i.e, the asymptotic limit of the entropy constructed with fundamental solution to the conjugate heat equation centered at $(x,t)$. All the definitions can be found in the next section. In this paper, we show that $\bar{W}(x,t)$ is a quantity that plays an equally important role as the asymptotic reduced volume in the study of ancient solutions.
\\

In Yokota \cite{yokota2009perelman}, he proved an important gap theorem for the asymptotic reduced volume of ancient solutions to the Ricci flow: the asymptotic reduced volume of every ancient solution that is not a Gaussian shrinker cannot be as large as some number strictly less than 1 and depending only on the dimension. Where the Gaussian shrinker is the Euclidean space evolving by the diffeomorphism generated by the gradient field of the function $\frac{|x|^2}{4}$ as well as rescaling.  Here we provide an entropy counterpart of Yokota's result, that is, a gap theorem for the asymptotic entropy on ancient solutions: if the asymptotic entropy of an ancient noncollapsed solution with bounded curvature is as large as some negative number, then the ancient solution must be a Gaussian shrinker. Our assumptions are stronger than Yokota's; this is because estimates for heat kernels are more difficult to obtain than estimates for the reduced distance. Our argument, although more complicated, resembles Yokota's in spirit. Moreover, in the special case of type \textmd{I} ancient solutions, according to Xu \cite{xuequation} it holds that $\bar{W}(x,t)=\log{\bar{V}(x,t)}$, where $\bar{V}(x,t)$ is the asymptotic reduced volume centered at $(x,t)$, then our gap theorem reduces to Yokota's.
\\

\begin{theorem_gap} \label{thm:main_2}
There exists $\varepsilon>0$ depending only on the dimension $n$ such that the following holds. Let $(M^n,g(t))_{t\in(-\infty,0]}$ be a complete ancient noncollapsed solution to the Ricci flow such that $\displaystyle \sup_{(x,t)\in M\times(-\infty,0]}|Rm|(x,t)<\infty$. If there exists $(x,t)\in M\times(-\infty,0]$ such that $\bar{W}(x,t)\geq-\varepsilon$, then $(M,g(t))$ is a Gaussian shrinker.\\
\end{theorem_gap}

In Perelman \cite{perelman2002entropy}, he asserted that on ancient solutions with nonnegative curvature operator, the condition that $\bar{W}(x,t)>-\infty$ is equivalent to $\kappa$-noncollapsing on all scales. This fact is widely known for the asymptotic reduced volume, but no similar result is known so far for the asymptotic entropy. In the following theorem, we show that under some moderate assumptions for an ancient solution with nonnegative curvature operator, bounded pointed entropy centered at some space-time point indeed implies noncollapsing on all scales at every previous time. This partially confirms Perelman's assertion. The implication in the other direction is proved in \cite{zhang2017on}. Notice in our argument, we do not assume $M$ is compact.
\\

\begin{theorem_noncollapsing} \label{thm:main_1}
Let $(M^n,g(t))_{t\in(-\infty,0]}$ be a complete ancient solution to the Ricci flow with nonnegative curvature operator. Assume $\displaystyle \sup_{x\in M}R(x,0)<\infty$, $\displaystyle \inf_{x\in M}\operatorname{Vol}(B_{g(0)}(x,1))>0$, where $R$ is the scalar curvature. Let $\bar{W}$ be the asymptotic entropy defined in (\ref{eq:W_bar_def}). If there exists a point $(x_0,t_0)$ such that $\bar{W}(x_0,t_0)\geq-\beta$, for some $\beta<\infty$, then there exists a $\kappa>0$, such that $(M,g(t))_{t\in(-\infty,t_0]}$ is $\kappa$-noncollapsed on all scales, where $\kappa$ depends only on $\beta$ and the dimension $n$.\\
\end{theorem_noncollapsing}

In this paper, unless otherwise specified, we always assume the Ricci flow that we work with is complete and has bounded geometry on every compact time interval, that is, the curvature is bounded from above and the volume of unit balls is bounded from below, notice the second assumption actually can be implied by the first and a noncollapsing assumption. For most of the time, we work with ancient solutions $(M,g(t))_{t\in(-\infty,0]}$ to the Ricci flow, on which it is known that the scalar curvature is nonnegative due to Chen \cite{chen2009strong}.
\\

This paper is organized as follows. In section 2, we provide some well known facts and fundamental tools that are necessary in our proof. In section 3, we show another version of Perelman's noncollapsing theorem, which is well suited to our argument. In section 4, we develop important estimates for the asymptotic Nash entropy; one may see the asymptotic Nash entropy has good dependence on its base point. In section 5, we use results from previous sections to complete the proof.

\bigskip

\section{Preliminaries}

\newtheorem{gradient_conjugate}{Lemma}[section]
\newtheorem{entropy_facts}[gradient_conjugate]{Lemma}
\newtheorem{log_sobolev}[gradient_conjugate]{Lemma}
\newtheorem{gradient_estimate}[gradient_conjugate]{Lemma}
\newtheorem{distance_distortion}[gradient_conjugate]{Lemma}

In this section, we will fix our notations and collect some known results that are very important to our proof. Let
\begin{eqnarray} \label{def:heat_kernel}
H_{(z,T)}(x,t)=\frac{1}{(4\pi(t-T))^{\frac{n}{2}}}\exp{(-f_{(z,T)}(x,t))},\ \ \text{where}\ T<t
\end{eqnarray}
be the fundamental solution to the heat equation coupled with the Ricci flow $(M^n,g(t))$, where $(z,T)$ is the base point. We may also fix $(x,t)$ and regard $H_{(z,T)}(x,t)$ as a function of $(z,T)$, then $H$ becomes the fundamental solution to the conjugate heat equation. The \textit{pointed entropy} and \textit{pointed Nash entropy} at time $T$ are defined as

\begin{eqnarray}
W_{(x,t)}(T)&=&\int_{M}[(t-T)(|\nabla_z f|^2+R(z,T))+f-n]Hdg_T(z), \label{eq:W_def}
\\
N_{(x,t)}(T)&=&\int_{M}fHdg_T(z)-\frac{n}{2},  \label{eq:N_def}
\end{eqnarray}
where $f=f_{(z,T)}(x,t)$, $H=H_{(z,T)}(x,t)$, and $dg_T(z)$ is the volume form for $g(T)$. The \textit{asymptotic entropy} and \textit{asymptotic Nash entropy} are defined as follows

\begin{eqnarray}
\bar{W}(x,t)&=&\lim_{T\rightarrow -\infty}W_{(x,t)}(T), \label{eq:W_bar_def}
\\
\bar{N}(x,t)&=&\lim_{T\rightarrow -\infty}N_{(x,t)}(T). \label{eq:N_bar_def}
\end{eqnarray}
where by monotonicity (see Lemma \ref{lem:basic} below), unless being negative infinity, these two limits always exist.
\\

We define some notations that are used in our argument. Let

\begin{eqnarray}
M_{T_1,T_2}=\sup_{x,y\in M}H_{(x,T_1)}(y,T_2),\ \ \text{where}\ T_1<T_2\leq 0. \label{eq:def_of_M}
\end{eqnarray}
We have that $M_{T_1,T_2}$ is finite, and by maximum principle,

\begin{eqnarray*}
M_{T_1',T_2'}\leq M_{T_1,T_2},\ \ \text{if}\ T_1'\leq T_1\ \text{and}\ T_2'\geq T_2.
\end{eqnarray*}
Notice that when applying the maximum principle to the conjugate heat equation, one needs to use the fact that $R\geq 0$ on an ancient solution. For $s<t$ and $x,y\in M$, define

\begin{eqnarray*}
d\nu_{(x,t)}^s(y)=H_{(y,s)}(x,t)dg_s(y).
\end{eqnarray*}
Clearly this is a probability measure since $\int_{M}H_{(y,s)}(x,t)dg_s(y)$ is independent of $s$ and $\displaystyle \lim_{s\rightarrow t+}H_{(y,s)}(x,t)=\delta_x(y)$.
\\

The following Hamilton-type \cite{hamilton1993matrix} gradient estimate for the conjugate heat equation shall provide justification for our work below.

\begin{gradient_conjugate}
Let $u$ be a positive solution to the conjugate heat equation coupled with the backward Ricci flow on $M^n\times [0,T]$, that is,
\begin{eqnarray*}
\frac{\partial}{\partial\tau}g&=&2Ric,
\\
\frac{\partial}{\partial\tau}u&=&\Delta u-Ru.
\end{eqnarray*}
Assume that the curvature and its first derivatives are bounded, and that $0<u\leq J$ on $M^n\times[0,T]$. Then there exists a constant $C$ depending only on the curvature and the curvature derivative upper bounds on $M\times[0,T]$, $T$ and dimension $n$ such that
\begin{eqnarray}
\tau\frac{|\nabla u|^2}{u^2}\leq C(1+\log{\frac{J}{u}})^2,
\end{eqnarray}
for all $\tau\in[0,T]$
\end{gradient_conjugate}

\begin{proof}
This lemma follows from Theorem 10 of \cite{ecker2008local}, where they have established that under the assumptions we make, the following holds
\begin{eqnarray*}
\frac{|\nabla u|^2}{u^2}(x,\tau)\leq(1+\log{\frac{J}{u}})^2(\frac{1}{\tau}+C+\frac{C\rho\coth(C\rho)+C}{\rho^2}),
\end{eqnarray*}
for all $x\in B_\tau(x_0,\rho)$ and $\tau\in(0,T]$, where $C$ is a constant depending on the curvature and the curvature derivative upper bounds and $x_0$ is a fixed point on $M$. The conclusion follows easily by taking $\rho$ to infinity.
\end{proof}

Under our assumption, that every Ricci flow we work with has bounded geometry on every compact time interval, we may apply the previous lemma to $H_{(z,T)}(x,t)$ on the interval $[T,\frac{T+t}{2}]$ and combine it with a Gaussian lower bound of the heat kernel (c.f. Theorem 26.31 of \cite{chow2010ricci}) to conclude that
\begin{eqnarray}
|\nabla_z f_{(z,T)}(x,t)|\leq C+C dist_{(g(T))}(z,x), \label{eq:growth_gradf}
\end{eqnarray}
where the constant $C$ depends on the curvature and the curvature derivative upper bounds and on the unit ball volume lower bound on the interval $[T,t]$, as well as on $t-T$. Moreover, it follows from the Bishop-Gromov comparison theorem that the volume growth is at most exponential. Therefore, we can justify integrability and integration by parts at infinity for all the formulae in Lemma \ref{lem:basic} below involving entropy that we shall work on. For instance, we can easily check that
\begin{eqnarray*}
\int_M \Delta_z f_{(z,T)}(x,t)H_{(z,T)}(x,t)dg_T(z)=\int_M |\nabla_z f_{(z,T)}(x,t)|^2H_{(z,T)}(x,t)dg_T(z).
\end{eqnarray*}
To see this, we consider the cut-off function $\displaystyle \psi(z)=\phi\left(\frac{d_{g(T)}(x,z)}{r}\right)$ where $\phi=1$ on $[0,1]$, $\phi$ decays smoothly to $0$ on $[1,2]$ with $|\phi'|<2$, and $\phi=0$ on $[2,\infty)$. It follows from integration by parts that
\begin{eqnarray*}
&&\int_M\psi(z)\Delta_z f_{(z,T)}(x,t)H_{(z,T)}(x,t)dg_T(z)
\\
&&=\int_M \psi(z)|\nabla_z f_{(z,T)}(x,t)|^2H_{(z,T)}(x,t)dg_T(z)+\int_M \langle\nabla_z\psi(z),\nabla_z f_{(z,T)}(x,t)\rangle H_{(z,T)}(x,t)dg_T(z),
\end{eqnarray*}
where
\begin{eqnarray*}
&&\left|\int_M \langle\nabla_z\psi(z),\nabla_z f_{(z,T)}(x,t)\rangle H_{(z,T)}(x,t)dg_T(z)\right|
\\
&&\leq \frac{2}{r}\int_{B_T(x,2r)\setminus B_T(x,r)}|\nabla_z f_{(z,T)}(x,t)|H_{(z,T)}(x,t)dg_T(z).
\end{eqnarray*}
Taking $r$ tending to infinity, the last term goes to zero because of the exponential volume growth, exponential quadratic decay of the heat kernel, as well as (\ref{eq:growth_gradf}); whence follows the integration by parts at infinity.
\\

We collect some well-known facts in the following lemma. One may find these results in Hein-Naber \cite{hein2014new}.

\begin{entropy_facts} \label{lem:basic}
The following hold for $T<t$:
\begin{enumerate}[(1)]
\item
  $\displaystyle \lim_{T\rightarrow t}W_{(x,t)}(T)=0$
\item
  $N_{(x,t)}(T)=\frac{1}{t-T}\int_T^tW_{(x,t)}(s)ds$
\item
  $W_{(x,t)}(T)\leq N_{(x,t)}(T)\leq0$
\item
  $\frac{d}{dT}W_{x,t}(T)=2(t-T)\int_{M}|Ric+\nabla^2f-\frac{g}{2(t-T)}|^2Hdg\geq 0$ (Perelman)
\item
  $\frac{d}{dT}N_{(x,t)}(T)=\frac{1}{t-T}(N_{(x,t)}(T)-W_{(x,t)}(T))\geq 0$
\item
  $N_{(x,t)}(T)=-\int_T^t 2(t-s)(1-\frac{t-s}{t-T})\int_{M}|Ric(z,s)+\nabla^2_z f_{(z,s)}(x,t)-\frac{g(z,s)}{2(t-s)}|^2d\nu_{(x,t)}^s(z)ds$
\item
  If there exists a $T<t$ such that $W_{(x,t)}(T)=0$, then the Ricci flow must be a Gaussian shrinker.\\
\end{enumerate}
\end{entropy_facts}

The following logarithmic Sobolev inequality and gradient estimate are very important for our arguments.

\begin{log_sobolev}[Hein-Naber \cite{hein2014new}]
Let $(M,g(t))_{t\in[0,T]}$ be a complete Ricci flow with bounded geometry on every compact time interval. For any $u\in C_0^{\infty}(M)$, the following hold:
\begin{eqnarray}
\int_{M}u^2d\nu-(\int_{M}ud\nu)^2&\leq&2(t-s)\int_{M}|\nabla u|^2d\nu, \label{eq:Poincare}
\\
\int_{M}u\log{u}d\nu-(\int_{M}ud\nu)\log{(\int_{M}ud\nu)}&\leq&(t-s)\int_{M}\frac{|\nabla u|^2}{u}d\nu, \label{eq:log_Sobolev}
\end{eqnarray}
where $0\leq s< t\leq T$, $d\nu=d\nu_{(x,t)}^s(y)$, and for the second inequality we also require $u\geq 0$.
\end{log_sobolev}

\textbf{Remark.} We will apply (\ref{eq:Poincare}) to the function $f$ as in definition (\ref{def:heat_kernel}), and apply (\ref{eq:log_Sobolev}) to the heat kernel. Even if the above inequalities are only proved for compactly supported functions, since at every time slice the heat kernel drops rapidly at infinity by our assumption on the Ricci flow, a cutoff argument easily justifies our work.
\\

The following gradient estimate was first proved by Hamilton \cite{hamilton1993matrix} on static manifolds; later it was adapted to the Ricci flow on closed manifolds by Zhang \cite{zhang2006some} and localized in \cite{bailesteanu2010gradient}. We will combine these results to obtain a lemma that can be used in our setting.
\\
\begin{gradient_estimate}[\cite{hamilton1993harnack}, \cite{zhang2006some}, \cite{bailesteanu2010gradient}] \label{lem:gradient_estimate}
Let $u:M^n\times [0,T]\rightarrow\mathbb{R}^+$ be a positive solution to the heat equation coupled with the Ricci flow, that is,
\begin{eqnarray*}
\frac{\partial}{\partial t}g&=&-2Ric,
\\
\frac{\partial}{\partial t}u&=&\Delta u.
\end{eqnarray*}
Assume the Riemann curvature is uniformly bounded and $0<u\leq J$ on $M^n\times[0,T]$. Then
\begin{eqnarray}
\frac{|\nabla u|^2}{u}\leq\frac{1}{t}u\log\frac{J}{u}. \label{eq:gradient}
\end{eqnarray}
\end{gradient_estimate}

\begin{proof}
Let $C_0$ be the bound of the Ricci curvature, that is, $|Ric|\leq C_0$ on $M\times [0,T]$. Applying Theorem 2.2 of \cite{bailesteanu2010gradient}, we have
\begin{eqnarray*}
\frac{|\nabla u|}{u}(x,t)\leq B(\frac{1}{2\rho}+\frac{1}{\sqrt{t}}+\sqrt{C_0})(1+\log\frac{J}{u}),\ \ \text{for all}\ x\in B_t(O,\rho)\ \text{and}\ t\neq 0,
\end{eqnarray*}
where $B$ is a constant depending only on the dimension of $M$ and $O$ is a fixed point on $M$.
\\

Taking $\rho\rightarrow\infty$, we have
\begin{eqnarray*}
t\frac{|\nabla u|^2}{u}\leq B^2(1+\sqrt{C_0t})^2(u^{\frac{1}{2}}+u^{\frac{1}{2}}\log\frac{J}{u})^2\leq C,\ \text{on}\ M\times[0,T],
\end{eqnarray*}
where $C$ is a constant depending on $B$, $C_0$, $T$ and $J$, and in the last inequality, we have also used the fact that the function $-x^{\frac{1}{2}}\log{x}$ is bounded from above and below on the interval $(0,J]$. Therefore, the subsolution to the heat equation (see \cite{zhang2006some})
\begin{eqnarray*}
t\frac{|\nabla u|^2}{u}-u\log\frac{J}{u}
\end{eqnarray*}
is bounded on $M\times[0,T]$, is nonpositive at time $0$, and we may apply the parabolic weak maximum principle (c.f. Theorem 12.10 of \cite{chow2008ricci}) to obtain the conclusion.
\end{proof}

The following distance distortion estimate is also helpful.

\begin{distance_distortion}[Perelman \cite{perelman2002entropy}]\label{lem:distant_distortion}
Suppose $Ric(x,t_0)\leq (n-1)K$ for all $x\in B_{t_0}(x_0,r_0)\bigcup\\ B_{t_0}(x_1,r_0)$. Then
\begin{eqnarray}
\frac{d}{dt}dist_t(x_0,x_1)\geq-2(n-1)(\frac{2}{3}Kr_0+r_0^{-1})\ \ \text{at}\ t=t_0. \label{eq:distance_distortion} \\\nonumber
\end{eqnarray}
\end{distance_distortion}

\bigskip

\section{Fundamental noncollapsing theorems}

\newtheorem{asymp_nash}{Lemma}[section]
\newtheorem{entropy_on_torus}[asymp_nash]{Corollary}
\newtheorem{prop_noncollapsing}[asymp_nash]{Proposition}
\newtheorem{perelman_noncollapsing}[asymp_nash]{Corollary}

We begin this section by computing the Nash entropy on $\mathbb{T}^k$, any flat torus of dimension $k\geq 1$. By Ricci-flatness, the Ricci flow on $\mathbb{T}^k$ is static and the Nash entropy of the Ricci flow coincides that of the linear heat equation. We will show that the Nash entropy cannot stay bounded; this fact has some consequence on the noncollapsing property which is important for our future development. It turns out that any collapsing sequence of parabolic cubes, after rescaling, will necessarily converge to a quotient space of $\mathbb{T}^k\times \mathbb{R}^{n-k}$, where $k\geq 1$, whose Nash entropy cannot stay bounded. Therefore, a lower bound of the Nash entropy implies a lower volume bound. In fact, all the properties we need here for $\mathbb{T}^k$ are the Ricci-flatness and the finiteness of the volume. \\

Let
\begin{eqnarray*}
H(x,y,t)=\frac{1}{(4\pi t)^{\frac{k}{2}}}e^{-f(x,y,t)}
\end{eqnarray*}
be the heat kernel on $\mathbb{T}^k$; $H$ coincides with the conjugate heat kernel coupled with the static backward Ricci flow. Then the Nash entropy centered at $(x,0)$ is defined as
\begin{eqnarray}
N(t)&=&\int_{\mathbb{T}^k}f(x,y,t)H(x,y,t)dy - \frac{k}{2} \label{def_static_nash}
\\\nonumber
&=&-\int_{\mathbb{T}^k}H(x,y,t)\log{H(x,y,t)}dy - \frac{k}{2}\log{(4\pi t)}-\frac{k}{2}.
\end{eqnarray}

\begin{asymp_nash} \label{lem:entropy_S}
$\displaystyle \lim_{t\rightarrow\infty}N(t)=-\infty$.
\end{asymp_nash}

\begin{proof}
From $H>0$ and $-H\log{H}\leq\frac{1}{e}$ for all $t>0$, it follows that
\begin{eqnarray*}
N(t)&=&-\int_{\mathbb{T}^k}H(x,y,t)\log{H(x,y,t)}dy - \frac{k}{2}\log{(4\pi t)}-\frac{k}{2}
\\
&\leq&\frac{1}{e}\operatorname{Vol}(\mathbb{T}^k)-\frac{k}{2}\log{(4\pi t)}-\frac{k}{2},\ \text{for all}\ t>0.
\end{eqnarray*}
The right-hand side obviously approaches negative infinity as $t$ approaches infinity.
\end{proof}

We have the following corollary.

\begin{entropy_on_torus} \label{coro:torus}
Let $N_{\mathbb{T}^k\times\mathbb{R}^{n-k}}(t)$ be the Nash entropy of the linear heat equation on the flat manifold $\mathbb{T}^k\times\mathbb{R}^{n-k}$, where $k\geq 1$. Then $\displaystyle \lim_{t\rightarrow\infty}N_{\mathbb{T}^k\times\mathbb{R}^{n-k}}(t)=-\infty$.
\end{entropy_on_torus}

\begin{proof}
It is easy to see that, if
\begin{eqnarray*}
H_1(x,t)&=&\frac{1}{(4\pi t)^\frac{n_1}{2}}\exp{(-f_1(x,t))},
\\
H_2(y,t)&=&\frac{1}{(4\pi t)^\frac{n_2}{2}}\exp{(-f_2(y,t))}
\end{eqnarray*}
are the heat kernels of Riemannian manifolds $(M_1,g_1)$, $(M_2,g_2)$ centered at $(x_0,0)$, $(y_0,0)$, respectively, then $H_1(x,t)H_2(y,t)$ and $N_{(x_0,0)}(t) + N_{(y_0,0)}(t)$ are the heat kernel and Nash entropy of $(M_1\times M_2,g_1+g_2)$ centered at $((x_0,y_0),0)$, respectively. It follows from Lemma \ref{lem:entropy_S} that $\displaystyle \lim_{t\rightarrow \infty} N_{\mathbb{T}^k\times\mathbb{R}^{n-k}}(t)=-\infty$.
\end{proof}

\begin{prop_noncollapsing} \label{prop:noncollapsing}
For all $\beta>0$, there exists $\kappa = \kappa (\beta,n)>0$ such that the following holds. Suppose a complete Ricci flow $(M^n,g(t))_{t\in[0,T]}$, $(x_0,t_0)\in M\times(0,T]$, and $r_0\in(0,\sqrt{t_0}\ ]$ are such that $|Rm|(x,t)\leq r_0^{-2}$ on $B_{t_0}(x_0,r_0)\times[t_0-r_0^2,t_0]$. If $N_{(x_0,t_0)}(t_0-r_0^2)\geq - \beta$, then $\operatorname{Vol}(B_{t_0}(x_0,r_0))\geq\kappa r_0^n$.
\end{prop_noncollapsing}

\begin{proof}
The statement is invariant under parabolic rescaling, so we may assume $r_0=1$. Suppose the conclusion does not hold. Then we can find a sequence of counterexamples whose injectivity radii at $(x_0,t_0)$ converge to $0$, while their curvatures on the parabolic cubes $B_{t_0}(x_0,1)\times[t_0-1,t_0]$ are bounded from above by $1$ and the Nash entropies satisfy $N_{(x_0,t_0)}(t)\geq -\beta$ for all $t\in[t_0-1,t_0)$. Further rescale this sequence of parabolic cubes to $\{B_{t_k}(x_k,r_k)\times[t_k-r_k^2,t_k]\}_{k=1}^\infty$ so that the injectivity radii at $(x_k,t_k)$ all become $1$, which implies that the curvature bounds will go to $0$ and the ``size'' $r_k$ of the parabolic cubes will go to infinity. Hence by \cite{hamilton1995compactness}, we can extract a subsequence converging in the pointed smooth Cheeger-Gromov sense to a flat ancient Ricci flow $(M^n_\infty,g_\infty(t),(x_\infty,t_\infty))_{t\in(-\infty,t_\infty]}$ whose injectivity radius at $(x_{\infty},t_{\infty})$ is $1$. It follows from Bieberbach's theorem (c.f. Theorem 98 of \cite{berger2012panoramic}) that $(M_\infty,g_\infty)\cong (\mathbb{T}^k\times\mathbb{R}^{n-k})/\Gamma$, where $k\geq 1$, $\mathbb{T}^k$ is a flat torus, and $\Gamma$ is a finite group.
\\
\newtheorem{claim_1}{Claim}

\begin{claim_1}
$N_{(x_\infty,t_\infty)}(T)\geq-\beta$, for all $T<t_\infty$.
\end{claim_1}

\begin{proof}[Proof of Claim 1]
For notational simplicity, we let $t_\infty=0$. Fix a $T<0$ and let $\{\Omega_i\times[T_i,0]\}_{i=1}^\infty$ be a compact space-time exhaustion of the limit flow. We assume without loss of generality that $T_i<T$ for any $i$. Let $g_i(t)$ and $H_i(z,t)=\frac{1}{(4\pi(-t))^\frac{n}{2}}\exp{(-f_i(z,t))}$ be the pull-back metrics and conjugate heat kernels, respectively, via diffeomorphisms given by the pointed Cheeger-Gromov convergence on $\Omega_i\times[T_i,0]$, where all the heat kernels are centered at the base points. According to \cite{lu2012convergence}, by passing to a subsequence the conjugate heat kernels centered at the base points will also converge to the heat kernel on the limit manifold centered at $(x_\infty,0)$, that is, $f_i\rightarrow f_\infty$ and $H_i\rightarrow H_\infty$ locally smoothly. By the assumption we made on the Nash entropy and Lemma \ref{lem:basic}(6), we have
\begin{eqnarray*}
\int_{T}^02s(1-\frac{s}{T})\int_{\Omega_i}|Ric_i(z,s)+\nabla_{g_i(z,s)}^2f_i(z,s)+\frac{g_i}{2s}|^2H_i(z,s)dg_i(z,s)ds\geq-\beta,\ \text{for all}\ i\in\mathbb{N}.
\end{eqnarray*}
Notice here we have used the fact that the integrands are nonpositive; hence their contributions to the integration outside $\Omega_i\times[-T_i,0]$ are also nonpositive. Passing this inequality to the limit and using Lemma \ref{lem:basic}(6) again for $(M_\infty,g_\infty(t))$ completes the proof.
\end{proof}

To apply this result, we need another claim.
\newtheorem{claim_2}[claim_1]{Claim}
\begin{claim_2}
Let $\pi:(\tilde{M},\tilde{g})\rightarrow (M=\tilde{M}/\Gamma,g)$ be a finite Riemannian cover, so that $|\Gamma|<\infty$. Then the Nash entropies defined as (\ref{def_static_nash}) of the linear heat equation on the cover and base satisfy
\begin{eqnarray*}
\tilde{N}_{(\tilde{x},0)}(t)\geq N_{(\pi(\tilde{x}),0)}(t).
\end{eqnarray*}
\end{claim_2}

\begin{proof}[Proof of Claim 2]
It is well known that
\begin{eqnarray*}
H(\pi(\tilde{x}),\pi(\tilde{y}),t)=\sum_{h\in\Gamma}{\tilde{H}(\tilde{x},h(\tilde{y}),t)}.
\end{eqnarray*}
We let $\tilde{M}_1\in\tilde{M}$ be a fixed fundamental domain of the covering $\tilde{M}$. It follows that
\begin{eqnarray*}
N_{(\pi(\tilde{x}),0)}(t) &=&-\int_{\tilde{M}_1}\sum_{h\in\Gamma}{\tilde{H}(\tilde{x},h(\tilde{y}),t)}\log{\left(\sum_{k\in\Gamma}{\tilde{H}(\tilde{x},k(\tilde{y}),t)}\right)}d\tilde{g}(\tilde{y})
-\frac{n}{2}\log{(4\pi t)}-\frac{n}{2}
\\
&\leq&-\int_{\tilde{M}_1}\sum_{h\in\Gamma}{\tilde{H}(\tilde{x},h(\tilde{y}),t)}\log{({\tilde{H}(\tilde{x},h(\tilde{y}),t)})}d\tilde{g}(\tilde{y})
-\frac{n}{2}\log{(4\pi t)}-\frac{n}{2}
\\
&=&-\int_{\tilde{M}}H(\tilde{x},\tilde{y},t)\log{H(\tilde{x},\tilde{y},t)}d\tilde{g}(\tilde{y})-\frac{n}{2}\log{(4\pi t)}-\frac{n}{2}
\\
&=&\tilde{N}_{(\tilde{x},0)}(t).
\end{eqnarray*}
\end{proof}

From this claim we deduce that the Nash entropy of $\mathbb{T}^k\times\mathbb{R}^{n-k}$ is bounded from below by $-\beta$, contradicting Corollary \ref{coro:torus}, thus completes the proof of Proposition \ref{prop:noncollapsing}\\
\end{proof}

\textbf{Remark.} Claim 2 in the previous proposition works also for solutions to the Ricci flow in place of static Riemannian manifolds. Note that we have only applied it to the static Ricci flat case.
\\

This proposition has a well-known consequence.
\\
\begin{perelman_noncollapsing}[Perelman's noncollapsing theorem]
Let $(M^n,g(t))_{t\in[0,T)}$ be a complete Ricci flow with $\nu[g(0),T]>-\infty$. Then there exists a positive number $\kappa$ depending only on $n$ and $\nu[g(0),T]$ such that the following holds. Let $(x_0,t_0)\in M\times (0,T)$ and $0<r_0\leq\sqrt{t_0}$ be such that $|Rm|(x,t)\leq r_0^{-2}$ for all $(x,t)\in B_{t_0}(x_0,r_0)\times[t_0-r_0^2,t_0]$. Then $\operatorname{Vol}(B_{t_0}(x_0,r_0))\geq\kappa r_0^n$.
\end{perelman_noncollapsing}

\begin{proof}
By items (3) and (5) of Lemma \ref{lem:basic}, we have
\begin{eqnarray*}
N_{(x_0,t_0)}(t_0-r_0^2)\geq N_{(x_0,t_0)}(0)\geq W_{(x_0,t_0)}(0)\geq\mu(g(0),t_0)\geq\nu[g(0),T].
\end{eqnarray*}
The conclusion follows from Proposition \ref{prop:noncollapsing}
\end{proof}

\bigskip

\section{Estimates on the asymptotic Nash entropy}

\newtheorem{decaying_bound}{Lemma}[section]
\newtheorem{gradient_converge}[decaying_bound]{Lemma}
\newtheorem{nash_invariant}[decaying_bound]{Proposition}
\newtheorem{mean_value_ineq}[decaying_bound]{Proposition}
\newtheorem{asymp_entropy}[decaying_bound]{Corollary}

In this section, we prove that under slightly weaker assumptions than those of the main theorems, the asymptotic Nash entropy $\bar{N}(x,t)$ not only is independent of the choice of $x$, but also decreases in $t$. This fact will be important for the proof of the main theorems. The estimates corresponding to the two main theorems, although different, are actually parallel to each other. Hence we will group these arguments together in the following development. Moreover, since every estimate constant depends on the dimension, we will suppress this dependance in our argument for simplicity.
\\
\begin{decaying_bound} \label{lem:decaying_bounf}
Let $(M^n,g(t))_{t\in(-\infty,0]}$ be a complete ancient solution to the Ricci flow.
\begin{enumerate}[(a)]
\item
  If there exists $C_0>0$ such that $\displaystyle\sup_{M\times(-\infty,0]}{|Rm|}\leq C_0$ and $\displaystyle\inf_{(x,t)\in M\times (-\infty,0]}\operatorname{Vol}(B_t(x,1))\geq C_0^{-1}$, then there exists a constant $C>0$ depending only on $C_0$ such that
  \begin{eqnarray*}
  0\geq N_{(x,t)}(T)\geq -\frac{n}{2}\log{(4\pi(t-T))}-C,
  \end{eqnarray*}
  for all $(x,t)\in M\times (-\infty,0]$ and $T<t-1$.
\item
  If $M$ has nonnegative curvature operator and if there exists $C_0>0$ such that $\displaystyle\sup_{x\in M}R(x,0)\\ \leq C_0$ $\displaystyle\inf_{x\in M}\operatorname{Vol}(B_0(x,1))\geq C_0^{-1}$ then there exists a constant $C>0$ depending only on $C_0$ such that
  \begin{eqnarray*}
  0\geq N_{(x,t)}(T)\geq -\frac{n}{2}\log{(4\pi(t-T))}-n\log{(-T+C)}-C,
  \end{eqnarray*}
  for all $(x,t)\in M\times (-\infty,0]$ and $T<t-1$.
\end{enumerate}
\end{decaying_bound}

\begin{proof}
In our argument we let $C$ be a constant depending only on $C_0$ and which may vary from line to line.
\begin{enumerate}[(a)]
  \item
  By our assumptions, we may apply a rough upper bound for the fundamental solutions of heat-type equations with respect to an evolving metric (c.f. Lemma 26.17 of \cite{chow2010ricci}) to conclude that
  \begin{eqnarray} \label{eq:bound_on_M_1}
  M_{T,T+1}\leq C(C_0),\ \text{for all}\ T\leq -1,
  \end{eqnarray}
  where $M_{T,T+1}$ is defined by (\ref{eq:def_of_M}). Note that the constant on the right-hand side does not depend on $T$.
  \\

  By the maximum principle applied to the heat equation and the conjugate heat equation, we have
  \begin{eqnarray*}
  H_{(x,T)}(y,t)\leq M_{T,T+1},\ \text{whenever}\ t\geq T+1,
  \\
  H_{(x,t)}(y,T)\leq M_{T-1,T},\ \text{whenever}\ t\leq T-1.
  \end{eqnarray*}
  Notice here that when applying the maximum principle to the conjugate heat equation, one should use the property that on an ancient solution to the Ricci flow the scalar curvature is always nonnegative. It then follows that
  \begin{eqnarray*}
  -\frac{n}{2}\log{(4\pi(t-T))}-f_{(z,T)}(x,t)=\log{H_{(z,T)}(x,t)}\leq \log{M_{T,T+1}}.
  \end{eqnarray*}
  Integrating this inequality against $d\nu_{(x,t)}^T(z)=H_{(z,T)}(x,t)dg_T(z)$, we have
  \begin{eqnarray}
  0\geq N_{(x,t)}(T)&\geq&-\frac{n}{2}-\log{M_{T,T+1}}-\frac{n}{2}\log{(4\pi(t-T))} \label{eq:nash_decay_a}
  \\\nonumber
  &\geq&-C-\frac{n}{2}\log{(4\pi(t-T))}.
  \end{eqnarray}
  Whence follows the conclusion.

  \item
  The proof is only a slight modification of part $(a)$. We need to estimate $M_{T,T+1}$ for this case. According to Hamilton's Harnack \cite{hamilton1993harnack}, we have $\displaystyle\sup_{M\times(-\infty,0]}|Rm|\leq C_0$. Applying Lemma \ref{lem:distant_distortion} with $K=C_0$ and $r_0=\frac{1}{2}$, we obtain that there exists $C>0$ depending only on $C_0$ such that $B_T(x,1-CT)\supseteq B_0(x,1)$, for all $T<0$ and $x\in M$. Furthermore, the Ricci flow with nonnegative scalar curvature shrinks the measure. Thus
  \begin{eqnarray*}
  \operatorname{Vol}(B_T(x,1-CT))\geq \operatorname{Vol}(B_0(x,1))\geq C_0^{-1}
  \end{eqnarray*}
  and
  \begin{eqnarray*}
  \operatorname{Vol}(B_T(x,1))\geq \frac{1}{C(1-CT)^n}
  \end{eqnarray*}
  by the Bishop-Gromov volume comparison theorem.
  Applying the rough upper bound for fundamental solutions of heat-type equations again (c.f. Lemma 26.17 of \cite{chow2010ricci}), we have
  \begin{eqnarray} \label{eq:bound_on_M_2}
  M_{T,T+1}\leq C(1-CT)^n.
  \end{eqnarray}
  Inserting this inequality into (\ref{eq:nash_decay_a}) completes the proof.\\
\end{enumerate}
\end{proof}

The following lemma is inspired by Hein-Naber's argument (see Theorem 4.5 of \cite{hein2014new}).
\begin{gradient_converge} \label{lem:gradient_converge}
Let $(M^n,g(t))_{t\in(-\infty,0]}$ be a complete ancient solution to the Ricci flow such that either:
\begin{enumerate}[(a)]
  \item
  There exists $C_0>0$ such that $\displaystyle\sup_{M\times(-\infty,0]}{|Rm|}\leq C_0$ and $\displaystyle\inf_{(x,t)\in M\times (-\infty,0]}\operatorname{Vol}(B_t(x,1))\geq C_0^{-1}$, or
  \item
  $M$ has nonnegative curvature operator and there exists $C_0>0$ such that $\displaystyle\sup_{x\in M}R(x,0)\leq C_0$ and $\displaystyle\inf_{x\in M}\operatorname{Vol}(B_0(x,1))\geq C_0^{-1}$.
\end{enumerate}
Then
\begin{eqnarray*}
\lim_{T\rightarrow-\infty}|\nabla_xN_{(x,t)}(T)|_{g(t)}=0
\end{eqnarray*}
for each $(x,t)\in M\times(-\infty,0]$. Moreover, the convergence is uniform in $x$.
\end{gradient_converge}

\begin{proof}
In the proof we always assume that $T<t-1$. By definition,
\begin{eqnarray*}
\nabla_xN_{(x,t)}(T)=\int_{M}[\nabla_xf_{(z,T)}(x,t)-f_{(z,T)}(x,t)\nabla_xf_{(z,T)}(x,t)]d\nu_{(x,t)}^T(z)
\end{eqnarray*}
and hence
\begin{eqnarray*}
|\nabla_xN_{(x,t)}(T)|_{g(t)}\leq \|\nabla_xf_{(z,T)}(x,t)\|_2(1+\|f_{(z,T)}(x,t)\|_2),
\end{eqnarray*}
by the Cauchy-Schwarz inequality, where the $L^2$ norms are taken with respect to $d\nu_{(x,t)}^T(z)$.\\

By Lemma \ref{lem:gradient_estimate} we have
\begin{eqnarray*}
|\nabla_xf_{(z,T)}(x,t)|_{g(t)}^2&=&\frac{|\nabla_xH_{(z,T)}(x,t)|_{g(t)}^2}{H_{(z,T)}(x,t)^2}
\\
&\leq&\frac{1}{t-T-1}\log{\frac{M_{T,T+1}}{H_{(z,T)}(x,t)}}
\\
&=&\frac{1}{t-T-1}(\log{M_{T,T+1}}+\frac{n}{2}\log{(4\pi(t-T))}+f_{(z,T)}(x,t)).
\end{eqnarray*}
Integrating against the measure $d\nu_{(x,t)}^T(z)$, we have
\begin{eqnarray} \label{eq:nashgrad_1}
\|\nabla_xf_{(z,T)}(x,t)\|_2^2\leq \frac{1}{t-T-1}(\log{M_{T,T+1}}+\frac{n}{2}\log{(4\pi(t-T))}+\frac{n}{2}),
\end{eqnarray}
where the last term on the right-hand side is because $N_{(x,t)}(T)\leq 0$.
\\

To estimate $\|f_{(z,T)}(x,t)\|_2$, we apply (\ref{eq:Poincare}) to obtain
\begin{eqnarray}
\|f_{(z,T)}(x,t)\|_2^2&=& \int_{M}f_{(z,T)}(x,t)^2 d\nu_{(x,t)}^T(z) \label{eq:nashgrad_2}
\\\nonumber
&\leq& \left(\int_{M}f_{(z,T)}(x,t)d\nu_{(x,t)}^T(z)\right)^2 + 2(t-T)\int_{M}|\nabla_zf_{(z,T)}(x,t)|_{g(T)}^2d\nu_{(x,t)}^T(z)
\\\nonumber
&\leq&(N_{(x,t)}(T)+\frac{n}{2})^2+n,
\end{eqnarray}
where the last inequality follows from
\begin{eqnarray*}
(t-T)\int_{M}(|\nabla_zf_{(z,T)}(x,t)|_{g(T)}^2+R(z,T))d\nu_{(x,t)}^T(z)\leq \frac{n}{2},
\end{eqnarray*}
which in turn is implied by the fact that $W_{(x,t)}(T)\leq N_{(x,t)}(T)$.
\\

Combining (\ref{eq:nashgrad_1}) and (\ref{eq:nashgrad_2}), we have
\begin{eqnarray} \label{eq:nashgrad_3}
|\nabla_xN_{(x,t)}(T)|_{g(t)}^2&\leq&\frac{1}{t-T-1}((N_{(x,t)}(T)+\frac{n}{2})^2+n)
\\\nonumber
&&\times(\log{M_{T,T+1}}+\frac{n}{2}\log{(4\pi(t-T))}+\frac{n}{2}).
\\\nonumber
\end{eqnarray}

Applying Lemma \ref{lem:decaying_bounf} as well as (\ref{eq:bound_on_M_1}) or (\ref{eq:bound_on_M_2}) to (\ref{eq:nashgrad_3}), we have:
\\

\emph{Case (a)}.
\begin{eqnarray*}
|\nabla_xN_{(x,t)}(T)|_{g(t)}^2\leq\frac{1}{t-T-1}(\frac{n}{2}\log{(4\pi(t-T))}+C)^3,
\end{eqnarray*}
or
\\

\emph{Case (b)}.
\begin{eqnarray*}
|\nabla_xN_{(x,t)}(T)|_{g(t)}^2\leq\frac{1}{t-T-1}(\frac{n}{2}\log{(4\pi(t-T))}+C+n\log{(-T+C)})^3,
\end{eqnarray*}
where $C$ depends only on $C_0$. In either case, the right-hand side converges to $0$ as $T\rightarrow-\infty$ whence follows the conclusion.
\\
\end{proof}

\begin{nash_invariant} \label{prop:nash_invariant}
Let $(M^n,g(t))_{t\in(-\infty,0]}$ be a complete ancient solution to the Ricci flow such that either:
\begin{enumerate}[(a)]
  \item
  there exists $C_0>0$ such that $\displaystyle\sup_{M\times(-\infty,0]}{|Rm|}\leq C_0$ and $\displaystyle \inf_{(x,t)\in M\times (-\infty,0]}\operatorname{Vol}(B_t(x,1))\geq C_0^{-1}$, or
  \item
  $M$ has nonnegative curvature operator and there exists $C_0>0$ such that $\displaystyle\sup_{x\in M}R(x,0)\leq C_0$ and $\displaystyle \inf_{x\in M}\operatorname{Vol}(B_0(x,1))\geq C_0^{-1}$.
\end{enumerate}
Then $\bar{N}(x,t)$ is independent of $x$.
\end{nash_invariant}

\begin{proof}
It suffices to show that if there exists $x\in M$ such that $\bar{N}(x,t)>-\infty$, then for any $y\in M$, $\bar{N}(y,t)=\bar{N}(x,t)$. By Lemma \ref{lem:gradient_converge}, in either case, it holds that for all $\varepsilon >0$, there exists $T_0\ll t$ such that $|\nabla_yN_{(y,t)}(T)|_{g(t)}<\varepsilon$ for all $y\in M$ and $T<T_0$. Hence $|N_{(y,t)}(T)-N_{(x,t)}(T)|\leq\varepsilon dist_{g(t)}(x,y)$ for all $y\in M$ and $T<T_0$. Taking $T\rightarrow -\infty$, we have $|\bar{N}(y,t)-\bar{N}(x,t)|\leq\varepsilon dist_{g(t)}(x,y)$. The proposition follows from taking $\varepsilon\rightarrow 0$.\\
\end{proof}

We have the following surprising integral inequality, which yields the monotonicity of the asymptotic Nash entropy.

\begin{mean_value_ineq} \label{prop:mean_value_ineq}
Let $(M^n,g(t))_{t\in(-\infty,0]}$ be a complete ancient solution to the Ricci flow such that either:
\begin{enumerate}[(a)]
  \item
  there exists $C_0>0$ such that $\displaystyle\sup_{M\times(-\infty,0]}{|Rm|}\leq C_0$ and $\displaystyle \inf_{(x,t)\in M\times (-\infty,0]}\operatorname{Vol}(B_t(x,1))\geq C_0^{-1}$, or
  \item
  $M$ has nonnegative curvature operator and there exists $C_0>0$ such that $\displaystyle\sup_{x\in M}R(x,0)\leq C_0$ and $\displaystyle \inf_{x\in M}\operatorname{Vol}(B_0(x,1))\geq C_0^{-1}$.
\end{enumerate}
Then for all $x,z\in M$ and $s<t\leq 0$ it holds that
\begin{eqnarray}
\bar{N}(x,t)\leq \int_{M}\bar{N}(y,s)d\nu_{(x,t)}^s(y)=\bar{N}(z,s).
\end{eqnarray}
Hence $\bar{N}(x,t)$ is independent of $x$ and is monotonically decreasing in $t$.
\end{mean_value_ineq}

\begin{proof}
Let $y,z\in M$ and fix a $T<0$ such that $T+1\ll s<t\leq 0$. Applying (\ref{eq:log_Sobolev}) to
\begin{eqnarray*}
u(y)&=&H_{(z,T)}(y,s)>0,
\\
d\nu&=&d\nu_{(x,t)}^s(y)=H_{(y,s)}(x,t)dg_s(y),
\end{eqnarray*}
we have
\begin{eqnarray*}
&&-H_{(z,T)}(x,t)\log{H_{(z,T)}(x,t)}
\\
&=&-(\int_{M}H_{(z,T)}(y,s)d\nu_{(x,t)}^s(y))\log{(\int_{M}H_{(z,T)}(y,s)d\nu_{(x,t)}^s(y))}
\\
&\leq&-\int_{M}H_{(z,T)}(y,s)\log{(H_{(z,T)}(y,s))}d\nu_{(x,t)}^s(y)
\\
&&+(t-s)\int_{M}\frac{|\nabla_yH_{(z,T)}(y,s)|^2}{H_{(z,T)}(y,s)}d\nu_{(x,t)}^s(y),
\end{eqnarray*}
where we have used the fact that

\begin{displaymath}
\int_{M}H_{(z,T)}(y,s)d\nu_{(x,t)}^s(y)=\int_{M}H_{(z,T)}(y,s)H_{(y,s)}(x,t)dg_s(y)=H_{(z,T)}(x,t).
\end{displaymath}
Integrating against $dg_T(z)$ and using definition (\ref{eq:N_def}), we have

\begin{eqnarray} \label{eq:mean_value_1}
&&N_{(x,t)}(T)+\frac{n}{2}+\frac{n}{2}\log(4\pi(t-T))
\\\nonumber
&\leq&\int_{M}N_{(y,s)}(T)d\nu_{(x,t)}^s(y)+\frac{n}{2}+\frac{n}{2}\log(4\pi(s-T))
\\\nonumber
&&+(t-s)\int_{M}\int_{M}\frac{|\nabla_yH_{(z,T)}(y,s)|^2}{H_{(z,T)}(y,s)}d\nu_{(x,t)}^s(y)dg_T(z),
\end{eqnarray}
where the last term needs to be estimated. By Lemma \ref{lem:gradient_estimate},

\begin{eqnarray*}
\frac{|\nabla_yH_{(z,T)}(y,s)|^2}{H_{(z,T)}(y,s)}\leq\frac{1}{t-T-1}[\log{M_{T,T+1}}-\log{H_{(z,T)}(y,s)}]H_{(z,T)}(y,s).\\
\end{eqnarray*}

Inserting this inequality into (\ref{eq:mean_value_1}), we have

\begin{eqnarray} \label{eq:mean_value_2}
N_{(x,t)}(T)&\leq&\int_{M}N_{(y,s)}(T)d\nu_{(x,t)}^s(y)+\frac{n}{2}\log{\frac{s-T}{t-T}}+\frac{t-s}{t-T-1}\log{M_{T,T+1}}
\\\nonumber
&&-\frac{t-s}{t-T-1}\int_{M}\int_{M}H_{(z,T)}(y,s)\log{(H_{(z,T)}(y,s))}d\nu_{(x,t)}^s(y)dg_T(z)
\\\nonumber
&=&(1+\frac{t-s}{t-T-1})\int_{M}N_{(y,s)}(T)d\nu_{(x,t)}^s(y)+\frac{t-s}{t-T-1}\log{M_{T,T+1}}
\\\nonumber
&&+\frac{t-s}{t-T-1}(\frac{n}{2}\log(4\pi(s-T))+\frac{n}{2})+\frac{n}{2}\log{\frac{s-T}{t-T}}.
\end{eqnarray}

By (\ref{eq:bound_on_M_1}) and (\ref{eq:bound_on_M_2}) we may observe that

\begin{eqnarray*}
\lim_{T\rightarrow-\infty}\frac{t-s}{t-T-1}\log{M_{T,T+1}}=0
\end{eqnarray*}
in either case (a) or case (b). Therefore, by taking $T\rightarrow-\infty$ in formula (\ref{eq:mean_value_2}), we have

\begin{eqnarray*}
\bar{N}(x,t)\leq\lim_{T\rightarrow-\infty}\int_{M}N_{(y,s)}(T)d\nu_{(x,t)}^s(y).
\end{eqnarray*}
Notice that $\{N_{(y,s)}(T)\}_{T\in(-\infty,s-1]}$ is a family of nonpositive functions of $y$ that are monotonic in $T$; the conclusion follows from Fatou's lemma as well as Proposition \ref{prop:nash_invariant}.\\

\end{proof}

To conclude this section, we remark that Proposition \ref{prop:nash_invariant} holds also for the asymptotic entropy, simply because the asymptotic entropy coincides with the asymptotic Nash entropy.

\begin{asymp_entropy}
Let $(M,g(t))_{t\in(-\infty,0]}$ be a complete ancient solution to the Ricci flow with bounded geometry on every compact time interval. Then
\begin{eqnarray*}
\bar{W}(x,0)=\bar{N}(x,0),
\end{eqnarray*}
for all $x\in M$. In particular, Proposition \ref{prop:mean_value_ineq} holds for the asymptotic entropy $\bar{W}(x,t)$.
\end{asymp_entropy}
\begin{proof}
By Lemma \ref{lem:basic}(2)(3)(4), for any $\varepsilon>0$ and $T<0$, we have 
\begin{eqnarray*}
W_{(x,0)}(T)\leq N_{(x,0)}(T)&=&-\frac{1}{T}\int_{T}^0W_{(x,0)}(s)ds\leq-\frac{1}{T}\int_{T}^{\varepsilon T}W_{(x,0)}(s)ds
\\
&\leq&-\frac{1}{T}\int_{T}^{\varepsilon T}W_{(x,0)}(\varepsilon T)ds=(1-\varepsilon)W_{(x,0)}(\varepsilon T).
\end{eqnarray*}
The conclusion follows from first taking $T\rightarrow-\infty$ and then $\varepsilon\rightarrow 0$.
\end{proof}

\bigskip

\section{Proof of the gap and noncollapsing theorems}
\newtheorem{main_prop_1}{Proposition}[section]
\newtheorem{corollary_noncollapsing}[main_prop_1]{Corollary}
\newtheorem{main_prop_2}[main_prop_1]{Proposition}

In this section, we present two propositions which imply our main theorems immediately via Lemma \ref{lem:basic} that $N_{(x,t)}(T)\geq W_{(x,t)}(T)$ and $\bar{N}(x,t)\geq\bar{W}(x,t)$.

\begin{main_prop_1}
Let $(M^n,g(t))_{t\in(-\infty,0]}$ be a complete ancient solution to the Ricci flow with nonnegative curvature operator. Assume $\displaystyle\sup_{x\in M}R(x,0)<\infty$, $\displaystyle \inf_{x\in M}\operatorname{Vol}(B_{g(0)}(x,1))>0$, where $R$ is the scalar curvature. Let $\bar{N}$ be the asymptotic Nash entropy defined in (\ref{eq:N_bar_def}). If there exists a point $(x_0,t_0)$ such that $\bar{N}(x_0,t_0)\geq-\beta$ for some $\beta<\infty$, then there exists a $\kappa>0$ such that $(M,g(t))_{t\in(-\infty,t_0]}$ is $\kappa$-noncollapsed on all scales, where $\kappa$ depends only on $\beta$ and the dimension $n$.
\end{main_prop_1}

\begin{proof}
Under the assumption of this proposition, we may apply part (b) of Proposition \ref{prop:mean_value_ineq} to conclude that $\bar{N}(x,t)\geq-\beta$, for all $t\leq t_0,\ x\in M$. Then we may apply Proposition \ref{prop:noncollapsing} to these points. So there exists a $\kappa=\kappa(n,\beta)>0$ such that for any $r>0$, if $|Rm|\leq r^{-2}$ on $B_t(x,r)\times[t-r^2,t]$, then $\operatorname{Vol}(B_t(x,r))\geq\kappa r^n$. That is, $(M,g(t))_{t\in(-\infty,t_0]}$ is $\kappa$-noncollapsed on all scales.\\
\end{proof}

The next corollary follows easily.
\\

\begin{corollary_noncollapsing}\label{corollary_noncollapsing}
Let $(M^n,g(t))_{t\in(-\infty,0]}$ be a complete ancient solution to the Ricci flow with the following assumptions: $\displaystyle\sup_{(x,t)\in M\times(-\infty,0]}|Rm|<\infty$, $\displaystyle \inf_{(x,t)\in M\times(-\infty,0]}\operatorname{Vol}(B_{g(t)}(x,1))>0$. If there exists a point $(x_0,t_0)$ such that $\bar{N}(x_0,t_0)\geq-\beta>-\infty$, then there exists a $\kappa>0$ such that $(M,g(t))_{t\in(-\infty,t_0]}$ is $\kappa$-noncollapsed on all scales, where $\kappa$ depends only on $\beta$ and the dimension $n$.
\end{corollary_noncollapsing}

\begin{proof}
The proof is similar to that of the previous proposition, where now part (b) is replaced by part (a) in Proposition \ref{prop:mean_value_ineq}.
\end{proof}

To prove the next proposition, one may apply a similar technique as used by Yokota \cite{yokota2009perelman}. He implemented a point picking method on the ancient solution to construct a contradicting sequence, which implies an $\varepsilon$-regularity theorem; his gap theorem follows from that. With all the tools we have developed by far, we are able to prove this theorem in a slightly different---if not significantly easier---way.

\begin{main_prop_2}
There exists $\varepsilon>0$ depending only on the dimension $n$ such that the following holds. Let $(M^n,g(t))_{t\in(-\infty,0]}$ be a complete ancient noncollapsed solution to the Ricci flow such that $\displaystyle\sup_{(x,t)\in M\times(-\infty,0]}|Rm|(x,t)<\infty$. If there exists $(x,t)\in M\times(-\infty,0]$ such that $\bar{N}(x,t)\geq-\varepsilon$, then $(M,g(t))$ is a Gaussian shrinker.\\
\end{main_prop_2}

\begin{proof}
Suppose $(M,g(t))$ is flat but not Euclidean. \ Then the noncollapsing assumption implies maximum volume growth. By Peter Li (see Corollary 16.3 of \cite{li2012geometric} for instance), $M$ is a finite quotient of Euclidean space. However, any finite group action $\Gamma\times \mathbb{R}^n\rightarrow \mathbb{R}^n$ has a fixed point, which is a contradiction since $M$ is smooth. To see this, one may take any $id\neq \gamma\in\Gamma$ and $x\in\mathbb{R}^n$; then $x+\gamma{x}+...+\gamma^{|\gamma|-1}{x}$ is a fixed point, where $|\gamma|>1$ is the order of $\gamma$. So, henceforth we assume $(M,g(t))$ is non-flat.
\\

Suppose the theorem is not true. Then we may find a sequence $\{(M_k,g_k(t),(x_k,0))\\_{t\in(-\infty,0]}\}_{k=1}^\infty$, such that the $(M_k,g_k(t))$ are non-flat and noncollapsed on all scales, $\displaystyle\sup_{M_k\times(-\infty,0]}\\|Rm_k|<\infty$, and $\bar{N}(x_k,0)\geq-\frac{1}{k}$. By the noncollapsing assumption, $\displaystyle \inf_{(x,t)\in M_k\times(-\infty,0]}\\ \operatorname{Vol}(B_t(x,1))>0$, and it follows from Corollary \ref{corollary_noncollapsing} that the sequence of ancient solutions are $\kappa$-noncollapsed with respect to a universal $\kappa$.
\\

For each $k$, we pick $(\tilde{x}_k,\tilde{t}_k)\in M_k\times(-\infty,0]$ such that $|Rm_k|(\tilde{x}_k,\tilde{t}_k)\geq\frac{1}{2}\displaystyle \sup_{M_k\times(-\infty,0]}|Rm_k|$. By parabolic rescaling the ancient flows centered at $(\tilde{x}_k,\tilde{t}_k)$ by the factors $|Rm_k|(\tilde{x}_k,\tilde{t}_k)$ with time shifts of $\tilde{t}_k$ to $0$, we obtain another sequence $\{(\tilde{M}_k,\tilde{g}_k(t),(\tilde{x}_k,0))_{t\in(-\infty,0]}\}_{k=1}^\infty$ which is $\kappa$-noncollapsed with respect to a universal $\kappa$, satisfies $\displaystyle \sup_{\tilde{M}_k\times(-\infty,0]}|\tilde{Rm}_k|<2$, $|\tilde{Rm}_k|(\tilde{x}_k,0)=1$, and $\bar{N}(\tilde{x}_k,0)\geq-\frac{1}{k}$. Here we have used the invariance of the Nash entropy under parabolic rescaling as well as Proposition \ref{prop:mean_value_ineq} to conclude that $\bar{N}(\tilde{x}_k,0)\geq-\frac{1}{k}$.
\\

By \cite{hamilton1995compactness}, we can extract a subsequence of $\{(\tilde{M}_k,\tilde{g}_k(t),(\tilde{x}_k,0))_{t\in(-\infty,0]}\}_{k=1}^\infty$ converging in the pointed Cheeger-Gromov sense to an ancient $\kappa$-noncollapsed Ricci flow $(M_\infty,g_\infty(t),\\ (x_\infty,0))_{t\in (-\infty,0]}$ with $|Rm_\infty|(x_\infty,0)=1$, where by \cite{lu2012convergence} the conjugate heat kernels centered at $(\tilde{x}_k,0)$ converge to the conjugate heat kernel centered at $(x_\infty,0)$. By the same argument as shown in Proposition \ref{prop:noncollapsing}, $N_{(x_\infty,0)}(T)\equiv 0$. Moreover, since the Nash entropy is the time average of pointed entropies (see Lemma \ref{lem:basic}), $W_{(x_\infty,0)}(T)\equiv 0$. It follows that $(M_\infty,g_\infty(t))$ is a Gaussian shrinker, contradicting $|Rm_\infty|(x_\infty,0)=1$.\\
\end{proof}

\textbf{Acknowledgement.}
The author would like to express his gratitude towards his doctoral advisors, Professor Bennett Chow and Professor Lei Ni, from whose vast knowledge in Ricci flow and patient mentorship he has continually benefited. He would also like to thank Professor Peng Lu for priceless suggestions and discussions.

\bibliographystyle{plain}
\bibliography{citation}

\noindent Department of Mathematics, University of California, San Diego, CA, 92093
\\ E-mail address: \verb"yoz020@ucsd.edu"

\end{document}